\documentclass[11pt]{article}
\usepackage{amsmath, epsfig, cite,colordvi,xcolor}
\usepackage{amssymb}
\usepackage{amsfonts}
\usepackage{latexsym}
\usepackage{graphicx}
\usepackage{tikz}
\usepackage{float}
\usepackage{mathrsfs}[12pt]
\usepackage{bm, amscd,  amsfonts, amssymb, cite,texdraw}
\usepackage{amsmath, epsfig, cite, color, colordvi}
\usepackage[colorlinks,urlcolor=black,linkcolor=black,anchorcolor=black,citecolor=black]{hyperref}
\usepackage{tikz}

\newtheorem{thm}{Theorem}

\newtheorem{lem}[thm]{Lemma}

\newcommand{\qed}{{\hfill\rule{4pt}{7pt}}}

\numberwithin{equation}{section}

\makeatletter \@addtoreset{equation}{section} \makeatother
\tikzstyle{every node}=[circle,inner sep=1pt,fill=white!60]
\tikzstyle{tn}=[shape=circle, draw, color=black!70]
\everymath{\displaystyle}
\begin{document}

\begin{center}
{\bf \large Some Identities Related to the Second-Order Eulerian Numbers }\\
\vspace{0.3cm}
Amy. M. Fu\\
School of Mathematics, Shanghai University of Finance and Economics\\
 Shanghai, 200433, China\\
Email: fu.mei@sufe.edu.cn
\end{center}

\noindent{\bf Abstract:}  We express the  N\"{o}rlund polynomials in terms of the second-order Eulerian numbers. Based on this expression, we derive several identities related to the Bernoulli numbers.  In particular, we present a short proof of the problem raised by Rz\c{a}dkowski and Urli\'{n}ska. 


\noindent{\bf AMS subject classfications:} 05A19; 05A15

\section{Introduction}
The Stirling permutations were introduced by Gessel and Stanley \cite{GesselandStanley}. For some related results on this subject, we refer to \cite{Bona,ChenandFu,MaYeh, JasonKubaandPahholzer}. Let $Q_n$ be the multiset $\{1,1,2,2,\ldots,n,n\}$.  A Stirling permutation of order $n$ is a permutation of $Q_n$ such that for each $1\leq m \leq n$, the elements lying between two occurrences of $m$ are greater than $m$. 

 The second-order Eulerian numbers $C_{n,k}$ count the Stirling permutations of order $n$ with $k$ decents, which satisfy the recurrence relation:
\begin{equation}\label{relation}
C_{n,k}=kC_{n-1,k}+(2n-k)C_{n-1,k-1}
\end{equation}
with $C_{1,1}=1$ and $C_{1,0}=0$.

By exhibiting a bijection between the set of partitions of $[n+m]$ with $m$ blocks and the bar permutations on the elements of $Q_n$ with $m$ bars, Gessel and Stanley \cite{GesselandStanley} proved that 
\begin{equation}\label{SC}
\sum_{m=0}^\infty S(n+m,m) x^m=\sum_{k=1}^n C_{n,k}x^k/(1-x)^{2n+1},
\end{equation}
where $S(n,m)$ are the Stirling numbers of the second kind.

The N\"{o}rlund polynomials $B^{(z)}_n$ can be defined by the exponential generating function:
\begin{equation}\label{Norlund}
\sum_{n=0}^\infty B^{(z)}_n \frac{x^n}{n!}=\left(\frac{x}{e^x-1}\right)^z.
\end{equation} 
Note that for fixed $n$, $B^{(z)}_n$ are polynomials in $z$ with degree $n$. If $z=0$, then we have $B^{(0)}_0=1$ and $B^{(0)}_{n}=0$ for $n\geq 1$. If $z=1$, then $B^{(1)}_n$ are the classical Bernoulli numbers $B_n$, i.e.,
\begin{equation}\label{Bernoulli}
    \sum_{n=0}^\infty B_n \frac{x^n}{n!}=\frac{x}{e^x-1}.
\end{equation}
If $z=n$, we have $B^{(n)}_n=(-1)^nc^{(2)}_n $, which are the Cauchy numbers 
of the second kind\cite{Comet}:
\begin{equation}\label{Cauchy}
    \sum_{n=0}^\infty c^{(2)}_n \frac{x^n}{n!}=\frac{-x}{(1-x)\ln (1-x)}.
\end{equation}

In \cite{Carlitz}, Carlitz showed that the  N\"{o}rlund polynomials and the Stirling numbers of the second kind satisfy the relation:
\begin{equation}\label{NorlundStirling}
S(m+n,m)={m+n \choose n} B^{(-m)}_{n}.
\end{equation}

Let $\langle z\rangle_n$ be the rising factorial defined by $\langle z\rangle_n=z(z+1)\cdots (z+n-1).$  The following result expresses the N\"{o}rlund polynomials in terms of the second-order Eulerian numbers.

\begin{thm}\label{theorem1} We have 
\begin{equation}\label{equation1}
    B^{(z)}_n=\frac{n!}{(2n)!}\sum_{k=1}^n (-1)^{k} C_{n,k}\langle z \rangle_k 
    \langle -z+n+1 \rangle_{n-k}.
\end{equation}
\end{thm}

Based on Theorem \ref{theorem1}, we obtain several identities involving the convolutions of Bernoulli numbers. In particular, we give a short proof of the problem raised by  Rz\c{a}dkowski and Urli\'{n}ska: 
\begin{equation}\label{RzUr}
\int_{0}^1 \sum_{k=0}^{n-1} C_{n,k+1} u^{k+1}(u-1)^{2n-k}du=\frac{B_{n+1}}{n+1}.
\end{equation}  
By computing the integral
$$
\int_{0}^1  u^{k+1}(u-1)^{2n-k}du=\frac{(-1)^k}{2(n+1)}{2n+1 \choose k}^{-1},
$$
we may restate \eqref{RzUr} as follows.
\begin{thm}\label{theorem3} We have 
    \begin{equation}\label{equation2}
        \sum_{k=1}^n (-1)^{k-1} {2n+1 \choose k}^{-1}C_{n,k}=2B_{n+1}.
    \end{equation}
\end{thm}
Notice that, by considering the partial derivative equation $\partial^n_t(\psi)=v_n(\psi)$, 
where $\psi$ is defined by the Lambert W-function $\psi(t,x)=W(xe^{x+t})$ and $v_n(x)$ is given by 
$$
v_n(x)=-x(1+x)^{-2n+1}\sum_{k\geq 1} (-1)^k C_{n,k}x^k,
$$  an alternative proof of \eqref{equation2} was discussed on Mathoverflow\cite{Mathflow} recently.

Let $H_n$ be the harmonic number defined by $H_n=\sum_{i=1}^n 1/i$. By using the $p$-adic method, Miki\cite{Miki} proved that 
\begin{equation}\label{Miki}
    \sum_{k=2}^{n-2} \frac{B_k}{k}\frac{B_{n-k}}{n-k}=2H_n\frac{B_n}{n}+\sum_{k=2}^{n-2}{n \choose k} \frac{B_k}{k}\frac{B_{n-k}}{n-k}.
\end{equation}
 
By computing the second derivative of the both sides of \eqref{equation1}, then employing Miki's identity \eqref{Miki}, we derive the following the result involving the harmonic numbers.
\begin{thm} \label{theorem4} We have 
    \begin{equation}\label{equation4}
    \sum_{k=1}^n (-1)^{k}  {2n-1 \choose k-1}^{-1}(H_{2n-k}-H_{k-1})C_{n,k}=
    \frac{n^2}{n-1}B_{n-1}+n\sum_{k=2}^{n-2} \frac{B_k}{k}\frac{B_{n-k}}{n-k}.
    \end{equation}
\end{thm}

 For $n \geq N$, Dilcher\cite{Dilcher} proved that 
\begin{equation}\label{Dilcher}
    \sum_{k_1+k_2+\cdots+k_N=n} {n \choose k_1, k_2,\ldots, k_N}B_{k_1}B_{k_2}\cdots B_{k_N}=N{n \choose N}\sum_{k=0}^{N-1}(-1)^{N-1-k}s(N,N-k)\frac{B_{n-k}}{n-k}, 
\end{equation}
where $s(n,k)$ is the Stirling number of the first kind. By setting $z=N$ in \eqref{equation1}, then employing  \eqref{Dilcher}, we obtain the following identity. 
\begin{thm}\label{theorem2} Let $N,n$ be nonnegative integers with $n \geq N$.  We have 
\begin{equation}\label{equation3}
\sum_{k=1}^n (-1)^k {2n-1 \choose N+k-1}^{-1}C_{n,k}=2n\sum_{k=0}^{N-1}(-1)^{N-1-k} s(N,N-k)\frac{B_{n-k}}{n-k}.
\end{equation}
In particular, if we let $n=N$ in \eqref{equation3}, the Cauchy numbers can be related by the second-order Eulerian numbers:
\begin{equation}\label{CauchyEulerian}
    2c^{(2)}_n=\sum_{k=1}^n (-1)^{n-k}{2n-1 \choose n+k-1}^{-1} C_{n,k}.
\end{equation}
\end{thm}

\section{Proofs}

If two polynomials in a single variable $z$ agree for every nonnegative integer $z$, then they agree as polynomials. Therefore, to prove Theorem \ref{theorem1}, it suffices to prove the following lemma.
\begin{lem}\label{lemma1} Given a nonnegative integer $m$, we have 
    \begin{equation}\label{NorE}
        B^{(-m)}_n=\frac{n!}{(2n)!}\sum_{k=1}^n(-1)^kC_{n,k}\langle-m\rangle_{k}\langle m+n+1\rangle_{n-k}.  
    \end{equation}
\end{lem}
\noindent{\it Proof:} Equating the coefficients of $x^m$ on the both sides of \eqref{SC} leads to 
\begin{equation}\label{EulerianStirling}
S(m+n,m)=\sum_{k=1}^{n}{2n+m-k \choose 2n}C_{n,k}.
\end{equation}
Comparing \eqref{NorlundStirling} and \eqref{EulerianStirling} for any integer $m \geq 0$, we have 
\begin{equation}\label{NES}
    B^{(-m)}_n=\frac{m!n!}{(m+n)!}\sum_{k=1}^{n}{2n+m-k \choose 2n}C_{n,k}=\frac{n!}{(2n)!}\sum_{k=1}^n(-1)^kC_{n,k}\langle-m\rangle_{k}\langle m+n+1\rangle_{n-k},    
\end{equation}
as desired. We complete the proof of Lemma \ref{lemma1}, and the proof of Theorem \ref{theorem1} as well.  \qed 

Observe that 
$$
\frac{d}{dx}\ln \left(\frac{e^x-1}{x}\right)=\frac{1}{x}\left(\frac{-x}{e^{-x}-1}-1\right)=\sum_{n=1}^\infty (-1)^n B_n \frac{x^{n-1}}{n!}.
$$
Thus,
$$
\ln \left(\frac{e^x-1}{x}\right)=\sum_{n=1}^\infty (-1)^n \frac{B_n}{n}\frac{x^n}{n!}=\frac{x}{2}+\sum_{n=2}^\infty \frac{B_n}{n}\frac{x^n}{n!}.
$$
Note that $B_1=-\frac{1}{2}$ and $B_n=0$ when $n$ is odd and greater than $1$. 

Since
$$
\frac{d^\ell}{dz^\ell}\left(\frac{x}{e^x-1}\right)^z
=\frac{d^\ell}{dz^\ell} e^{z\ln (\frac{x}{e^x-1})}=(-1)^\ell e^{z\ln (\frac{e^x-1}{x})}\left[\ln \left(\frac{e^x-1}{x}\right)\right]^\ell,
$$
we have 
\begin{equation}\label{Dequation1}
    \frac{d}{dz}B^{(z)}_{n}|_{z=0}=-\frac{B_n}{n}
\end{equation}
and 
\begin{equation}\label{Dequation2}
    \frac{d^2}{dz^2}B^{(z)}_{n}|_{z=0}=\frac{n}{n-1}B_{n-1}+\sum_{k=2}^{n-2}{n \choose k} \frac{B_k}{k}\frac{B_{n-k}}{n-k}.
\end{equation}
Using Miki's identity \eqref{Miki}, we can rewrite \eqref{Dequation2} as 
\begin{equation}\label{Dequation3}
    \frac{d^2}{dz^2}B^{(z)}_{n}|_{z=0}=\frac{n}{n-1}B_{n-1}+\sum_{k=2}^{n-2}\frac{B_k}{k}\frac{B_{n-k}}{n-k}-2H_n\frac{B_n}{n}.
\end{equation}

\vspace{0.3cm}

\noindent{\it Proof of Theorem \ref{theorem3}:} Applying the derivation $d/dz|_{z=0}$ to both sides of \eqref{equation1}, then using \eqref{Dequation1}, we have 
\begin{equation}\label{Mequation}
-\frac{B_n}{n}=\frac{1}{2n}\sum_{k=1}^n (-1)^{k}{2n-1 \choose k-1}^{-1} C_{n,k}.
\end{equation}
Mutiplying both sides by $-2n$, then replacing $n$ by $n+1$, we have 
\begin{equation}\label{MEq}
\sum_{k=1}^{n+1}(-1)^{k-1}{2n+1 \choose k-1}^{-1}C_{n+1,k}=2B_{n+1}.
\end{equation}
By the recurrence relation \eqref{relation}, 
the left-hand side of \eqref{MEq} can be rewritten as 
\begin{eqnarray*}
\lefteqn{\sum_{k=1}^{n+1}(-1)^{k-1}{2n+1 \choose k-1}^{-1}\left(kC_{n,k}+(2n+2-k)C_{n,k-1}\right)}\\
&&=\sum_{k=1}^{n}(-1)^{k-1}\frac{k!(2n+2-k)!}{(2n+1)!}C_{n,k}-\sum_{k=1}^{n}(-1)^{k-1}\frac{k!(2n+1-k)!(2n+1-k)}{(2n+1)!}C_{n,k}\\
&&=\sum_{k=1}^{n}(-1)^{k-1}\frac{k!(2n+1-k)!}{(2n+1)!}C_{n,k}=
\sum_{k=1}^n (-1)^{k-1}{2n+1 \choose k}^{-1}C_{n,k},
\end{eqnarray*}
as desired. We complete the proof of Theorem \ref{theorem3}. \qed

\noindent{\it Proof of Theorem \ref{theorem4}:} Apply the derivation $d^2/dz^2|_{z=0}$ to both sides of \eqref{equation1}. By \eqref{Mequation}, we have 
\begin{eqnarray*}
    \frac{d^2}{dz^2}B^{(z)}_{n}|_{z=0}
&=&\frac{1}{n} \sum_{k=1}^n (-1)^k {2n-1 \choose k-1}^{-1} C_{n,k}(H_{k-1}-(H_{2n-k}-H_n))\\
&=&\frac{1}{n} \sum_{k=1}^n (-1)^k {2n-1 \choose k-1}^{-1} C_{n,k}(H_{k-1}-H_{2n-k})-2H_n\frac{B_n}{n}.
\end{eqnarray*}
Combining with \eqref{Dequation3}, after some arrangements, we complete the proof of Theorem \ref{theorem4}. \qed

We conclude this paper by the {\it proof of Theorem \ref{theorem2}.}  Setting $z=N$ in \eqref{Norlund} leads to
\begin{eqnarray*}
\sum_{k=0}^\infty B^{(N)}_k\frac{x^k}{k!}&=&\left(\frac{x}{e^x-1}\right)^N=\left(\sum_{k=0}^\infty B_k\frac{x^k}{k!}\right)^N\\
&=&\sum_{n=0}^\infty \frac{x^n}{n!} \sum_{k_1+k_2+\cdots+k_N=n}{n \choose k_1,k_2,\ldots,k_N}B_{k_1}B_{k_2}\cdots B_{k_N},
\end{eqnarray*}
where $k_1,k_2,\ldots, k_N$ are nonnegative integers. 

Hence,
$$
B^{(N)}_n=\sum_{k_1+k_2+\cdots+k_N=n}{n \choose k_1,k_2,\ldots, k_N}B_{k_1}B_{k_2}\cdots B_{k_N}.
$$
By Theorem \ref{theorem1}, for $1\leq N\leq n$, we have 
\begin{eqnarray*}
    B^{(N)}_n&=& \frac{n!}{(2n)!} \sum_{k=1}^n (-1)^k C_{n,k} \langle N\rangle_k \langle-N+n+1\rangle_{n-k}\\
    &=&\frac{N}{2n}{n \choose N} \sum_{k=1}^n (-1)^k C_{n,k}{2n-1 \choose N+k-1}^{-1}.
\end{eqnarray*}
Combining  with Dilcher's identity \eqref{Dilcher}, we have
$$
\frac{N}{2n}{n \choose N} \sum_{k=1}^n (-1)^k C_{n,k}{2n-1 \choose N+k-1}^{-1}=N{n \choose N}\sum_{k=0}^{N-1} (-1)^{N-1-k} s(N,N-k)\frac{B_{n-k}}{n-k}.
$$
Divided both sides by $\frac{N}{2n}{n \choose N}$, we complete the proof of Theorem \ref{theorem2}. \qed

\end{document}